\documentclass[a4paper,11pt]{article}
\usepackage{MyLaTeX,MyArticle}
\usepackage[vcentermath,enableskew]{youngtab}
\usepackage{multirow}
\input xy
\xyoption{all}

\title{Lower Bounds for Representation Growth}
\author{David A.~Craven, University of Oxford}
\date{5th October, 2009}

\begin{document}
\maketitle

\begin{abstract} This article examines lower bounds for the representation growth of finitely generated (particularly profinite and pro-$p$) groups. It also considers the related question of understanding the maximal multiplicities of character degrees in finite groups, and in particular simple groups.
\end{abstract}

\section{Introduction}

The representation growth of residually finite (particularly profinite) groups is a relatively new area of research (\cite{lubotzkymartin2004}, \cite{larsenlubotzky2008}), but hopefully will result in a theory as rich as the study of subgroup growth \cite{lubotzkysegal}. In this work we study lower bounds for the representation growth of pro-$p$ and profinite groups, and the connected topic of character degree multiplicities for finite groups.

In \cite{craven2008}, the author proved a result on the character degrees of the symmetric groups that, when combined with work in \cite{moreto2007} and \cite{jaikin2003}, yields a corollary that might be thought of as a basic result in representation growth. Let $G$ be a finitely generated, residually finite group, and let $r_n(G)$ be the number of inequivalent, complex irreducible representations of $G$ of dimension $n$, whose kernels have finite index. The main aim of representation growth is to relate the arithmetic properties of the sequence $\(r_n(G)\)$ with algebraic properties of the group. It is well known \cite[Proposition 2]{blmm2002} that all of the $r_n(G)$ are finite if and only if all finite-index subgroups of $G$ have finite abelianizations (the so-called FAb property).

\begin{thma}\label{theorema} There exists a function $f:\N\to\N$ such that $f(x)\to\infty$ as $x\to\infty$ with the following property: for any infinite, finitely generated, residually finite group $G$ with FAb, the sequence $r_n(G)$ is greater than $f(n)$ infinitely often. In particular, for any such $G$, the sequence $\(r_n(G)\)$ contains a divergent subsequence.
\end{thma}

This has the rather satisfactory consequence that the representation growth of a group does indeed grow. (The corresponding statement for subgroup growth is, of course, false.) Theorem \ref{theorema} in particular proves that representation growth cannot be arbitrarily slow.
%

If only finitely many (but at least one) of the $r_n(G)$ are infinite then one possibility is that we could simply ignore the first finitely many of the $r_n(G)$ and study the remaining sequence. The next theorem deals with this situation.

\begin{thma}\label{theoremb} Let $G$ be an infinite, finitely generated, residually finite group. Let $I(G)$ be the set of all natural numbers $i$ such that $r_i(G)=\infty$. The following are equivalent:
\begin{enumerate}
\item $I(G)$ is finite and non-empty;
\item $r_i(G)=0$ for all sufficiently large $i$; and
\item $G$ is virtually abelian.
\end{enumerate}
\end{thma}

Combining these two results gives the following corollary, producing the four possible types of residually finite group that can exist.

\begin{cora}\label{corollaryc} Let $G$ be a finitely generated, residually finite group. Let $I(G)$ be the set of all natural numbers $i$ such that $r_i(G)=\infty$. Exactly one of the following possibilities holds:
\begin{enumerate}
\item $I(G)=\emptyset$ and only finitely many of the $r_i(G)$ are non-zero;
\item $I(G)=\emptyset$ and the sequence $(r_i(G))$ contains a subsequence that tends to infinity;
\item $0<|I(G)|<\infty$, and only finitely many of the $r_i(G)$ are non-zero; and
\item $I(G)$ is infinite.
\end{enumerate}
In the first case, $G$ is finite, in the second, $G$ has FAb and is infinite, in the third case $G$ is infinite and virtually abelian, and in the final case $G$ has none of the previous properties.
\end{cora}

Groups in the final category include (infinite) finitely generated nilpotent groups, which can be studied using so-called \emph{twist isoclasses} \cite[Theorem 6.6]{lubotzkymagid1985}, and so even in this class something might be said. However, using purely the $r_n(G)$ and their related objects (like partial sums, zeta functions, and so on) it is likely that nothing much can be said about such groups.

\bigskip

As mentioned earlier in the introduction, we said that we will study the representations of symmetric groups. The process in \cite{craven2008} to generate irreducible characters with the same degree is constructive, and in Section \ref{secsymgrps} we derive an explicit bound, proving the following result.

\begin{thma}\label{explicitsym} Let $n$ be an integer, and let $X(n)$ denote the multiset of the degrees of the irreducible characters of the symmetric group $S_n$. Let $m(n)$ denote the largest of the multiplicities of the elements of $X(n)$. 
For all sufficiently large $i$, if $n\geq 3\cdot 81^i$ then $m(n)\geq 2^i$. In particular, for all sufficiently large $n$,
\[ m(n)\geq n^{1/7}.\]
\end{thma}

In Section \ref{secsymgrps} we derive a much more complicated explicit bound, but $n^{1/7}$ is approximately correct; so (for large $n$) the symmetric group $S_n$ has of the order of (approximately) $n^{0.16}$ characters of the same degree. Of course, this is only a lower bound, but it is hard to believe that the real answer isn't $\bigO{n^\ep}$ for some $\ep$. 

\begin{conja}\label{symrepgrowth} Let $n$ be an integer, and let $X(n)$ denote the multiset of the degrees of the irreducible characters of the symmetric group $S_n$. Let $m(n)$ denote the largest of the multiplicities of the elements of $X(n)$. There are positive constants $\ep_1$ and $\ep_2$ such that, for all sufficiently large $n$,
\[ n^{\ep_1}<m(n)<n^{\ep_2}.\]
\end{conja}

Theorem \ref{explicitsym} provides the lower bound in Conjecture \ref{symrepgrowth}. Notice that, when one changes from the \emph{degree} of the symmetric group to the \emph{order} of the symmetric group, then all of the functions $n^\ep$ become functions of the same order. Thus Conjecture \ref{symrepgrowth} would give the exact order of the growth of the maximal multiplicity of symmetric group character degrees in terms of the order of the group. This agrees with \cite[Corollary 2.7]{liebeckshalev2004} (see also \cite[Corollary 1.4(ii)]{liebeckshalev2005}), which implies that the growth of the maximal multiplicity in terms of the \emph{order} is slower than $n^\ep$ for any $\ep>0$. In Section \ref{simplegroups}, we compare this growth to the other simple groups, and prove that the alternating groups are the finite simple groups $G$ for which the maximum $m(G)$ of the $r_n(G)$ grows slowest relative to the order of the group $G$.

In Section \ref{propgroups} we turn our attention to $p$-groups, using the known results on conjugacy classes of $p$-groups to derive bounds for the growths of $r_n(G)$ and the partial sums $R_n(G)$. After considering $r_n(G)$ and $R_n(G)$ for powerful pro-$p$ groups and the Nottingham group, we consider all pro-$p$ groups. The strongest result we derive here is an easy consequence of a remarkable recent theorem of Jaikin-Zapirain \cite{jaikin2009un} (stated here as Theorem \ref{jaikinamazing}), which resolves a problem first posed by Pyber in \cite{pyber1992}.

\begin{thma}\label{repgrowthpgroup} There exists a constant $c$ such that, if $G$ is a finitely generated pro-$p$ group, then for all sufficiently large $n$,
\[ R_{p^n}(G)\geq cn\frac{\log_p n}{\log_p \log_p n}.\]
and for infinitely many $n$,
\[ r_{p^n}(G)\geq 2c\frac{\log_p n}{\log_p \log_p n}.\]
\end{thma}

The original proof of the second part of the result above -- that $r_{p^n}(G)$ has a divergent subsequence -- first appeared in \cite{jaikin2003}, but the new results in \cite{jaikin2009un} produce a (much) better bound.

In the final section, we discuss the concept of a lower bound and provide examples, constructed by Kassabov and Nikolov in \cite{kassabovnikolov2006}, that show that the concept of a lower bound needs to be modified to match that of Theorem \ref{theorema}.

\begin{thma}\label{repgrowthminimum} Let $f$ be a divergent, non-decreasing function. There exists a finitely generated profinite group $G$ such that $R_n(G)<f(n)$ infinitely many times. In other words, there is no divergent, non-decreasing function $f$ such that $f(n)<R_n(G)$ for any finitely generated profinite group $G$, and for all sufficiently large $n$.
\end{thma}

Theorems \ref{theorema} and \ref{repgrowthminimum} taken together indicate that in general, the representation growth of finitely generated profinite groups can behave very badly, and so it might be necessary to consider only certain classes of profinite group, such as pro-$p$ groups, for which things seem much better.

\section{Proof of Theorem \ref{theorema}}

The main tool for the proof of Theorem \ref{theorema} is the following result from \cite{craven2008}, itself depending on results from \cite{jaikin2003} and \cite{moreto2007}.

\begin{thm}[{{\cite[Corollary 1.3]{craven2008}}}] \label{finitegroupthm}Let $G$ be a finite group, and let $m(G)$ denote the maximum of $r_n(G)$, where $n\in \N$. Then $|G|$ is bounded by a function of $m(G)$.
\end{thm}

We will briefly mention how this theorem is proved. It relies on two special cases of this result, for $p$-groups and for finite simple groups. For $p$-groups, this is Theorem \ref{repgrowthpgroup}, and for simple groups this is the results from Sections \ref{secsymgrps} and \ref{simplegroups}. Using these two explicit computations, we can firstly give a bound for soluble groups, and then extend this to all finite groups using the generalized Fitting subgroup.

Write $\ell:\N\to\N$ for the function
\[ \ell(n)=\min_{|G|\geq n} m(G).\]
By Theorem \ref{finitegroupthm}, this function is well defined and non-decreasing. Also, for every finite group $G$ of order $n$, we have that $r_i(G)\geq \ell(n)$ for some $i<\sqrt n$. In fact, since $|G|\geq r_j(G)j^2$ for all $j$, we see that
\[ i\leq \sqrt{\frac{n}{\ell(n)}}.\]

Thus if $H$ is a quotient of $G$, and $H$ has order $n$, then there exists $i\leq\sqrt{n/\ell(n)}$ such that $r_i(G)\geq \ell(n)$.

Now let $G$ be an infinite, finitely generated, residually finite group, and suppose that $G$ has the FAb property. There is an infinite sequence $H_1,H_2,\dots$ of finite quotients of $G$ with $|H_i|<|H_{i+1}|$ for all $i$. For $H_i$ of order $n_i$, let $j_i$ denote the natural number, at most $\sqrt{n_i/\ell(n_i)}$, such that $r_{j_i}(H_i)\geq \ell(n_i)$. Therefore
\[ r_{j_i}(G)\geq \ell(n_i),\]
and thus the $r_{j_i}(G)$ form a divergent subsequence of the $r_m(G)$, bounded below by $\bar\ell(n_i)$, proving Theorem 
\ref{theorema}.

\bigskip

This theorem proves that the sequence $(r_n(G))$ contains a divergent subsequence, but for the $R_n(G)$, we can get reasonable growth bounds that are probably not massively far from the truth. Using a theorem of Pyber from \cite{pyber1992}, it is easy to show that
\[ R_n(G)\geq c\frac{\log n}{(\log\log n)^8}\]
for infinitely many $n$.

However, a more involved argument, due to Andrei Jaikin-Zapirain, proves something considerably better.

\begin{prop}\label{propRnGprofinite} There is a constant $c>0$ such that if $G$ is a finitely generated profinite group with FAb, then there are infinitely many integers $n$ for which
\[ R_n(G)\geq c\log n(\log\log n)^{1-\ep},\]
for any $\ep>0$.
\end{prop}
\begin{pf} Let $G$ be a finitely generated profinite group. If $G$ possesses infinitely many maximal subgroups, then $G$ maps onto infinitely many finite groups with trivial Frattini subgroup; in particular, onto such groups with arbitrarily large order. At the end of the proof of Theorem A from \cite{pyber1992}, Pyber proves that if $H$ is a finite group with trivial Frattini subgroup, then
\[ k(H)\geq 2^{c (\log |H|)^{1/8}}\]
for some constant $c>0$; in particular, all the quotients $H$ with trivial Frattini subgroup satisfy $R_{|H|}(H) \geq 2^{c(\log |H|)^{1/8}}$, and so therefore does $G$ for infinitely many $n=|H|$. This is well above the bound needed.

Thus we may assume that $G$ possesses only finitely many maximal subgroups, in which case $G$ is virtually pro-nilpotent, with pro-nilpotent subgroup $H$ of finite index. In this case, Theorem \ref{repgrowthpgroup} proves that $R_n(H)$ grows at least as quickly as $c\log n(\log \log n)^{1-\ep}$ for any $\ep>0$, and hence so does $R_n(G)$, as required. (To move between $R_n(H)$ and $R_n(G)$, we note that if $N$ is a normal subgroup of a finite group $G$, then $k(G)\geq k(N)/|G:N|$ (see e.g., \cite[Lemma 2.1(ii)]{pyber1992}), so if $|G:H|$ is fixed, $R_n(G)$ grows with the same order as $R_n(H)$.)
\end{pf}

Notice that it seems that pro-nilpotent groups (or pro-$p$ groups) are the bounding case in this result. In particular, if $G$ is a finitely generated groups with lots of simple groups, then $2^{c(\log n)^{1/8}}$ is the slowest that $R_n(G)$ can grow, at least for infinitely many $n$. We will return to this concept of only being able to bound infinitely many $R_n(G)$ from below in Section \ref{alternatinggroups}.






\section{Proof of Theorem \ref{theoremb}}

We start with a lemma, which gives us extra information in the case where $r_1(G)=\infty$.

\begin{lem}\label{infmanylinearchars} Let $G$ be a finitely generated group. If $r_1(G)=\infty$ then for all $n$, either $r_n(G)=0$ or $r_n(G)=\infty$.
\end{lem}
\begin{pf} Suppose that $G$ has a representation $\phi$ of degree $n$, with kernel $K$, and let $\psi$ be a $1$-dimensional representation, with kernel $H$. If $\phi\otimes \psi=\phi$, then it must be that $H$ contains $K$, since else the kernel of $\phi\otimes\psi$ would not be $K$. Thus for each representation $\phi$ of degree $n$, there are only finitely many $1$-dimensional representations $\psi$ such that $\phi\otimes\psi=\phi$, and so $r_n(G)=\infty$, as claimed.
\end{pf}

We now prove Theorem \ref{theoremb}. Firstly, if $G$ is virtually abelian, then it has some infinite abelian subgroup $H$ of index $n$. If $\rho$ is a representation of $G$, then an irreducible constituent $\phi$ of $\rho\res H$ is $1$-dimensional, and so by Frobenius reciprocity $\rho$ has dimension at most $n$, giving (iii) implies (ii); that (ii) implies (i) is clear.

Now suppose that $I(G)$ is finite and non-empty; since $G$ is not FAb, choose a subgroup $H$ of index $n$ with $|H/H'|$ infinite.

Suppose that there are infinitely many $i$ such that $r_i(H)=\infty$, and write $I=I(H)$. Let $J=I(G)$ be the corresponding (finite) set for $G$. We will derive a contradiction, proving that $I$ is finite.

Let $\rho$ be an irreducible representation of $H$, of dimension $m$. Since $(\rho\ind G)\res H$ has $\rho$ as a constituent, there is some constituent $\psi$ of $\rho\ind G$ such that $\rho$ is a constituent of $\psi\res H$; thus $m\leq \dim \psi$. Note that $\rho\ind G$ has dimension $nm$, and so $\dim \psi\leq nm$.

Let $a\in I$ be greater than any element of $J$. There are infinitely many (inequivalent) representations $\rho_i$ of $H$ of dimension $a$, and thus there must be infinitely many representations $\psi_i$ such that $\rho_i$ is a constituent of $\psi_i\res H$ and $\dim \psi_i$ lies between $a$ and $am$. Since there are only finitely many constituents of a given representation, this implies that $r_c(G)$ is infinite for some $a\leq c\leq am$, a contradiction. Thus $I(H)$ is finite, and so by Lemma \ref{infmanylinearchars} $r_i(H)$ is non-zero for only finitely many $i$; let $c$ be the largest dimension of an irreducible representation of $H$. By a well-known theorem of Jordan \cite{jordan1878}, for each finite quotient $H/K$ of $H$, there is an abelian normal subgroup $W/K$ such that $|H/W|$ is bounded by $d=f(c)$ for some non-decreasing function $f$. Since $H$ is finitely generated, there are only finitely many subgroups of index at most $d$, so let $A$ be the intersection of all such subgroups, necessarily a normal subgroup of finite index in $H$. We claim that $A$ is residually (finite abelian), and is hence abelian. This proves that $H$, and hence $G$, is virtually abelian, as required.

Let
\[ H=H_1\geq H_2\geq H_3\geq \cdots\]
be a descending sequence of normal subgroups of finite index of $H$ such that $\bigcap H_i=1$. Let $A_i=A\cap H_i$, and note that the $A_i$ is a descending sequence of normal subgroups of finite index of $A$ such that $\bigcap A_i=1$. Since $H/A_i$ is a finite group, there is some abelian subgroup $B/A_i$ of index at most $d$, and by construction $B\geq A$, so $A/A_i$ is abelian, as required.

\bigskip

Corollary \ref{corollaryc} follows from Theorems \ref{theorema} and \ref{theoremb}, of course.

\section{Degree Multiplicity for Symmetric Groups}
\label{secsymgrps}
This section relies on work of the author in \cite{craven2008}, and we will briefly recall what is involved there. There is a standard bijection between the irreducible characters of the symmetric group of degree $n$ and the partitions of $n$, with the degree of a particular character calculable from the corresponding partition, via \emph{hook numbers}. We presume that the reader is familiar with this method, and we will pause to fix notation only.

If $\lambda$ is a partition, write $|\lambda|$ for the number of which $\lambda$ is a partition, and $\lambda'$ for the conjugate of $\lambda$. Let $t(\lambda)$ be the sum of the number of rows of $\lambda$ and the number of columns of it. Let $H(\lambda)$ denote the multiset of all hook numbers of $\lambda$. If $H(\lambda)=H(\mu)$, then the characters corresponding to $\lambda$ and $\mu$ have the same degree.

To any partition, one may associate the \emph{enveloping partition} $E(\lambda)$, which is constructed in \cite{craven2008}, and is illustrated here by example. One takes a square of length $t(\lambda)$, appends a copy of $\lambda$ to the left and to the bottom of the square, and removes a reflected copy of $\lambda$ from the bottom-right portion of the square. This is the enveloping partition of $(5,3,3,2)$:

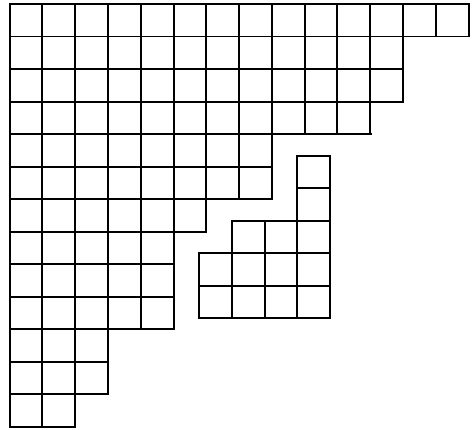
\[ \hspace{-0.6in}\yng(14,12,12,11,8,8,6,5,5,5,3,3,2)\hspace{-1.4in}\raisebox{-0.5in}{\young(:::\,,:::\,,:\,\,\,,\,\,\,\,,\,\,\,\,)}\]

\medskip

If $\lambda$ is a partition, and $t=t(\lambda)$, then write $E(\lambda)_i$ for the partition got from $E(\lambda)$ by incrementing the first $t$ rows of $E(\lambda)$ by $i$. In \cite{craven2008}, the following theorem is proved.

\begin{thm}[{{\cite[Theorem 7.1]{craven2008}}}] Suppose that $\lambda$ and $\mu$ are partitions, and that $H(\lambda)=H(\mu)$. Write $t$ for the sum of the number of rows and the number of columns of $\lambda$. (This is the same as that for $\mu$.) Then
\[ H(E(\lambda)_i)=H(E(\mu)_i).\]
\end{thm}

If we start with a partition $\lambda$ that is not self-conjugate, then the partitions $\lambda$ and $\lambda'\neq\lambda$ have the same hook numbers. From these two partitions, we may construct four partitions with the same hook numbers, namely
\[ E(\lambda)_1,\ E(\lambda')_1,\ E(\lambda)_1',\text{ and }E(\lambda')_1'.\]
If $|\lambda|=n$ and $t=t(\lambda)$, then all four of these partitions are partitions of $n+t^2+t$. This procedure can be iterated, to produce, given a partition $\lambda$, a set of $2^i$ partitions with the same hook numbers. Here we will calculate the smallest integer $N$ such that it can be guaranteed using this procedure that for all $n\geq N$, there are $2^i$ different partitions of $n$, each of which has the same hook numbers.

Firstly, given a partition $\lambda$ with $|\lambda|=n$ and $t(\lambda)=t$, we need to calculate the size of the partition got by applying the above procedure of taking $\mu\mapsto E(\mu)_1$ a number, say $i$, of times. Let $f$ denote the function on the set of all partitions given by $f(\mu)=E(\mu)_1$, and write $n_1=n$ and $t_1=t$. It is clear that
\[ |E(\lambda)_1|=n+t+t^2\text{ and }t(E(\lambda)_1)=3t+1.\]
Therefore, if $n_i$ and $t_i$ denote the size and row and column sum of $f^{(i-1)}(\lambda)$ (i.e., $f$ applied $i-1$ times to $\lambda$), we see that
\[ t_i=3t_{i-1}+1,\qquad n_i=n_{i-1}+t_{i-1}+t_{i-1}^2.\]
The first recurrence is easily solved to get
\[ t_i=3^{i-1}t+\frac{3^{i-1}-1}2,\]
and solving the second recurrence relation yields
\[ n_i=n+\frac{(4t^2+4t+1)(9^{i-1}-1)}{32}-\frac{i-1}{4}.\]

The equations above imply that given a partition $\lambda$ that is not self-conjugate with $|\lambda|=n$ and $t(\lambda)=t$, one may construct $2^i$ partitions with row and column sum $t_i$, and by extending the first $t_{i-1}$ rows by $j$ each, they may be taken to have sizes $n_i+jt_{i-1}$ for all $j\geq 0$.

The idea is to find $t_{i-1}$ partitions, each of which has the same row and column sum $t$, and whose sizes cover the $t_{i-1}$ congruence classes modulo $t_{i-1}$. Therefore for some integer $N$ we would have found $2^i$ partitions of size $N$ with the same hook numbers, and for all subsequent integers as well.

Suppose that a partition $\lambda$ has $t(\lambda)=t$. Furthermore, suppose that $t$ is odd (so that $\lambda$ is definitely not self-conjugate), and write $t=2r+1$. Then the largest that $|\lambda|$ can be is $(t^2-1)/4$ (which is a rectangle of sides $r$ and $r+1$), and the smallest that $|\lambda|$ can be is $t-1$ (which is a hook). Furthermore, it is easy to see that every possible size between these two can be given by a partition that is not self-conjugate. Thus given a row and column sum $t$, there are $(t^2-4t+7)/4$ different possibilities for $n$, and these possiblities form an interval.

Finally, since $t_{i-1}=3^{i-2}t+(3^{i-2}-1)/2$, we see that there are enough partitions if
\[ t^2-4t+7\geq 4\cdot 3^{i-2}t+2(3^{i-2}-1).\]
Using the quadratic formula, we get the exact solution
\[ t=2(1+3^{i-2})\pm \sqrt{4(1+3^{i-2})^2+2\cdot 3^{i-2}-9},\]
and we take the approximate solution
\[ t=5+4\cdot 3^{i-2},\]
which guarantees that there are enough partitions. Notice that the smallest value of $n$ is $t-1$, and therefore substituting these values into the equation for $n_i$ gives
\[ n_i=4+4\cdot 3^{i-2}+\frac{(4(5+4\cdot 3^{i-2})^2+4(5+4\cdot 3^{i-2})+1)(9^{i-1}-1)}{32}-\frac{i-1}{4}.\]

Thus we have proved that the symmetric group $S_n$ has $2^i$ irreducible characters of the same degree if 
\[n\geq \frac{15-16\cdot 3^{i-1}+1025\cdot 9^{i-2}+1584\cdot 27^{i-2}+576\cdot 81^{i-2}-8i}{32},\]
as required.

This is far from optimal. In \cite{craven2008}, it was shown that for all $n\geq 22$ there are four partitions with the same hook numbers, whereas this strategy proves it only for $n\geq 98$. For eight partitions, this method requires $n\geq 3078$, and while the real bound is not known precisely, it is known to be true for $n\geq 200$. In general, however, there appears to be no easy improvement on the method above.

\section{The Other Finite Simple Groups}
\label{simplegroups}

Apart from the alternating groups, the finite simple groups are the sixteen classes of groups of Lie type, together with the twenty-six sporadic simple groups. In terms of asymptotic group theory, the sporadic groups are unimportant, but we briefly mention the maximal degree multiplicities of the sporadic groups in a table, derived from the information in \cite{atlas}.

\begin{center}\begin{tabular}{|c|c||c|c||c|c|}
\hline $G$ & $m(G)$ & $G$ & $m(G)$ & $G$ & $m(G)$
\\ \hline $M_{11}$ & $3$ & $Co_3$ & $3$ &  $B$ & $2$
\\ $M_{12}$ & $3$ & $McL$ & $2$ & $M$ & $3$
\\ $M_{22}$ & $2$ & $Suz$ & $3$ & $J_1$ & $3$
\\ $M_{23}$ & $3$ & $He$ & $3$ & $ON$ & $3$
\\ $M_{24}$ & $3$ & $HN$ & $3$ & $J_3$ & $3$
\\ $HS$ & $3$ & $Th$ & $2$ & $Ru$ & $3$
\\ $J_2$ & $2$ & $Fi_{22}$ & $4$ & $J_4$ & $3$
\\ $Co_1$ & $2$ & $Fi_{23}$ & $3$ & $Ly$ & $5$
\\ $Co_2$ & $3$ & $Fi_{24}'$ & $2$ & $T$ & $2$
\\\hline
\end{tabular}\end{center}

[We include the Tits group $T={}^2F_4(2)'$ here, since it is `semi-sporadic', and not really one of the Ree groups $^2F_4(2^{2n+1})$.] What is interesting here is that, with the exception of the Lyons and smallest Fischer groups, all of the sporadic groups have maximal multiplicity either $2$ or $3$. In particular, if $G$ is a finite simple group and $m(G)=2$, then $|G|$ is at most that of the Baby Monster, and if $m(G)=3$, then $|G|$
is at most that of the Monster.

For the alternating groups, it is easy to see that
\[ \frac25 m(S_n)\leq m(A_n)\leq \frac52 m(S_n),\]
using Clifford theory.

\begin{lem}\label{cliffordmG} Let $G$ be a finite $p$-group, and let $N$ be a normal subgroup of index $p^n$. We have
\[ \(\frac{p}{p^2+1}\)^nm(G)\leq m(N)\leq \(\frac{p^2+1}{p}\)^nm(G).\]
\end{lem}
\begin{pf} Firstly assume that $n=1$. Suppose that there are $m$ characters of $G$ of the same degree. There are $pi$ of them that restrict to $i$ irreducible characters $\psi$ of $N$ (as $p$ of them restrict to each such $\psi$), and $m-pi$ of them that do not restrict to an irreducible character, and instead restrict to $p(m-pi)$ characters of $N$ (with the same degree).

The case where there are the fewest characters of the same degree in $N$ is when $i=p(m-pi)$, and so $i=mp/(p^2+1)$. Hence $m(N)\geq m(G)\cdot p/(p^2+1)$. However, by Frobenius reciprocity the situation is exactly the same for induction from $N$ to $G$, and so $m(G)\geq m(N)\cdot p/(p^2+1)$. A simple induction completes the proof.
\end{pf}

Of course, this leaves only the groups of Lie type, so fix a Lie-type group $G=G(q)$. It is known \cite[Theorem 1.7]{liebeckshalev2005} that the orders, character degrees, and character degree multiplicities, of $G$ are polynomials in $q$ (dependent on the type, but just the Lie rank of the group determines a lot). For the exceptional groups, these polynomials are known, and have been collated by Frank L\"ubeck; they are currently available on his website \cite{lubeckwebsite}. However, these data are only for the adjoint and simply connected versions of the groups, and so the simple group is not given for $^\ep E_6(q)$ and $E_7(q)$ (for certain values of $q$). However, using some elementary Clifford theory and the tables of character degrees, it is possible to still get the maximal multiplicities for the simple groups.

For the general group of Lie type, the polynomials are in the order $q$ of the finite field over which the group lies. For the Suzuki and Ree groups -- $^2B_2(q)$, $^2G_2(q)$, and $^2F_4(q)$ -- we use the notation $q^2=p^{2n+1}$, where $p$ is either $2$ or $3$. In \cite{liebeckshalev2005}, Liebeck and Shalev prove that if the Lie rank of $G(q)$ is $\ell$, then $m(G(q))$ is a polynomial with degree $\ell$; in the tables below, we reproduce the \emph{exact} polynomial for $m(G(q))$ (and the character degree at which it is attained) for each of the exceptional groups, and describe afterwards the small values of $q$ for which the table is incorrect.

For the \emph{simple} group $G(q)$, with $q$ odd, the table is as below.
\begin{footnotesize}\begin{center}\begin{tabular}{|c|c|c|}
\hline Group & Degree & Multiplicity
\\ \hline $G_2(q)$ & $q^6-1$ & $(q-1)^2/2$
\\ $^2G_2(q)$ & $q^6+1$ & $(q^2-3)/2$
\\ $^3D_4(q)$ & $(q^6-1)(q^4-q^2+1)(q^2-q+1)$ & $(q^4-2q+1)/4$
\\ $F_4(q)$ & $(q^{12}-1)(q^8-1)(q^2-1)(q^2-q+1)$ & $q^2(q^2-1)/6$
\\ $E_6(q)$ & $(q^{12}-1)(q^9-1)(q^6-1)(q^5-1)(q^4-1)$ & $(q^4-1)(q^2-1)/8\gcd(q-1,3)$
\\ $^2E_6(q)$ & $(q^{18}-1)(q^{12}-1)(q^{10}-1)(q^6-1)(q^4-1)/(q^9-1)(q^5-1)$ & $(q^4-1)(q^2-1)/8\gcd(q+1,3)$
\\ \multirow{2}{*}{$E_7(q)$} & $(q^{18}-1)(q^{12}-1)(q^{10}-1)(q^8-1)(q^7-1)(q^6-1)(q^2-1)$ & \multirow{2}{*}{$q(q^6-1)/28$}
\\ & $(q^{18}-1)(q^{14}-1)(q^{12}-1)(q^{10}-1)(q^8-1)(q^6-1)(q^2-1)/(q^7-1)$ & 
\\ $E_8(q)$ & $(q^{30}-1)(q^{24}-1)(q^{20}-1)(q^{18}-1)(q^{14}-1)(q^{12}-1)(q^2-1)$ & $(q^4-1)(5q^4-2q^3-7)/64$
\\ \hline
\end{tabular}\end{center}\end{footnotesize}
For $q$ even, it looks like this.
\begin{footnotesize}\begin{center}\begin{tabular}{|c|c|c|}
\hline Group & Degree & Multiplicity
\\ \hline $^2B_2(q)$ & $q^4+1$ & $(q^2-2)/2$
\\ $G_2(q)$ & $q^6-1$ & $q(q-2)/2$
\\ $^3D_4(q)$ & $(q^6-1)(q^4-q^2+1)(q^2-q+1)$ & $q(q^3-2)/4$
\\ $F_4(q)$ & $(q^{12}-1)(q^8-1)(q^2-1)(q^2-q+1)$ & $q^2(q^2-1)/6$
\\ $^2F_4(q)$ & $(q^{24}-1)(q^4+1)/(q^4+q^2+1)$ & $q^2(q^2-2)/4$
\\ $E_6(q)$ & $(q^{12}-1)(q^9-1)(q^6-1)(q^5-1)(q^4-1)$ & $q^4(q^2-1)/8\gcd(q-1,3)$
\\ $^2E_6(q)$ & $(q^{18}-1)(q^{12}-1)(q^{10}-1)(q^6-1)(q^4-1)/(q^9-1)(q^5-1)$ & $q^4(q^2-1)/8\gcd(q+1,3)$
\\ \multirow{2}{*}{$E_7(q)$} & $(q^{18}-1)(q^{12}-1)(q^{10}-1)(q^8-1)(q^7-1)(q^6-1)(q^2-1)$ & \multirow{2}{*}{$q(q^6-1)/14$}
\\ & $(q^{18}-1)(q^{14}-1)(q^{12}-1)(q^{10}-1)(q^8-1)(q^6-1)(q^2-1)/(q^7-1)$ & 
\\ $E_8(q)$ & $(q^{30}-1)(q^{24}-1)(q^{20}-1)(q^{18}-1)(q^{14}-1)(q^{12}-1)(q^2-1)$ & $q^4(5q^4-2q^3-8)/64$
\\ \hline
\end{tabular}\end{center}\end{footnotesize}

We should describe briefly how to determine the values in the table for $^\ep E_6(q)$ and $E_7(q)$ when there is a non-trivial centre. For $E_6(q)$ for example, there are $m=(q^4-1)(q^2-1)/8$ characters $\phi$ of the degree $f$ in the table for the adjoint group $E_6(q).3$, and since there are no characters of degree $(\deg\phi)/3$ for the simply connected group $3\cdot E_6(q)$, there are $m/3$ characters of degree $f$ for the simple group $E_6(q)$. Also, if $\psi$ is a character of $E_6(q).3$ then there are no characters of degree $\deg\psi/3$, and so it suffices to consider those character degrees for $E_6(q).3$ whose multiplicities exceed $m/3$, and in all cases it is easy to see that one gets fewer than $m$ characters with the same degree for the simple group. The technique is similar for $^2E_6(q)$ and $E_7(q)$.

Note that the maximal multiplicity of character degrees for $E_7(q)$ is realized by two different sets of characters, as suggested in the table: the first one has the smaller degree, and is also more naturally expressed as a product of polynomials of the form $(q^i-1)$.

There are obviously some small exceptions, and these are summarized below. The only unresolved case is $E_7(3)$, for which the multiplicity lies betwen $78$ and $80$. Na\"ive Clifford theory and the information for the adjoint and simply connected versions of $E_7(3)$ appears to be not enough to determine the multiplicities.

\begin{center}\begin{tabular}{|ccc|ccc|}
\hline Group & Degree & Multiplicity & Group & Degree & Multiplicity
\\ \hline $^2B_2(8)$ & $35$ & $3$ & $E_7(2)$ & $5070690584338804425$ & $9$
\\ $^3D_4(2)$ & $351$ & $3$ & $F_4(2)$ & $541450$ & $4$
\\ $E_6(2)$ & $42826799925$ & $8$ & $G_2(2)'$ & $7$ & $3$
\\ $E_6(3)$ & $127752132719411200$ & $84$ & $G_2(3)$ & $91$ & $3$
\\ $^2E_6(2)$ & $27498621150$ & $5$ & $G_2(4)$ & $819$ & $7$
\\ \hline
\end{tabular}\end{center}

For the classical groups, there is no known general formula for the maximal multiplicity of the character degrees, and so we use the lower bounds given in \cite{moreto2007}. [The choice of $d$ in the table below is influenced by the requirement that the numerator in each multiplicity should be the same.]

\begin{center}\begin{tabular}{|c|c|c|}
\hline Group & $\Orth(|G|)$ & Multiplicity
\\ \hline $\PSL_d(q)$ & $\frac{q^{d^2-1}}{\gcd(q-1,d)}$ & $\vphantom{\D\(\frac{\phi(q^n-1)}{n^2(q-1)}\)}\frac{\phi(q^d-1)}{d^2(q-1)}$
\\ $\PSU_{2d}(q)$ & $\frac{q^{4d^2-1}}{\gcd(q+1,2d)}$ & $\vphantom{\D\(\frac{\phi(q^n-1)}{n^2(q-1)}\)}\frac{\phi(q^d-1)}{4d^2}$
\\ $\PSU_{2d+1}(q)$ & $\frac{q^{4d(d+1)}}{\gcd(q+1,2d+1)}$  & $\vphantom{\D\(\frac{\phi(q^n-1)}{n^2(q-1)}\)}\frac{\phi(q^d-1)}{(2d+1)^2}$
\\ $\PSp_{2d}(q)$ & $\frac{q^{2d^2+d}}{\gcd(2,q-1)}$ & $\vphantom{\D\(\frac{\phi(q^n-1)}{n^2(q-1)}\)}\frac{\phi(q^d-1)}{4d}$
\\ $\POmega_{2d+1}(q)$ & $\frac{q^{2d^2+d}}{\gcd(2,q-1)}$ & $\vphantom{\D\(\frac{\phi(q^n-1)}{n^2(q-1)}\)}\frac{\phi(q^d-1)}{4d+2}$
\\ $\POmega^+_{2d}(q)$ & $\frac{q^{2d^2-d}}{\gcd(4,q^d-1)}$ & $\vphantom{\D\(\frac{\phi(q^n-1)}{n^2(q-1)}\)}\frac{\phi(q^d-1)}{4d}$
\\ $\POmega^-_{2d+2}(q)$ & $\frac{q^{2d^2+d+1}}{\gcd(4,q^{d+1}+1)}$ & $\vphantom{\D\(\frac{\phi(q^n-1)}{n^2(q-1)}\)}\frac{\phi(q^d-1)}{4d+4}$
\\ \hline\end{tabular}\end{center}

We aim to prove that each of these grows faster (in terms of $|G|$) than the symmetric groups can possibly do, proving that the symmetric groups are definitely the simple groups with the slowest-growing function $m(G)$ in terms of $|G|$. We will prove that, for sufficiently large $|G|$, for symmetric groups, 
\[ \log(\log m(G)+\log\log |G|)<(\log\log |G|)/2,\]
whereas for groups of Lie type the opposite inequality holds. This shows that the growth of $m(G)$ with respect to $|G|$ is slower for the symmetric groups than for the groups of Lie type, proving our claim.

We firstly note that the number of partitions of $m$ is asymptotically
\[ p(m)\sim \frac{\e^{a\sqrt m}}{bm},\]
where $a=\pi\sqrt{2/3}$ and $b=4\sqrt3$, by the famous Hardy--Ramanujan asymptotic formula \cite[(1.41)]{hardyramanujan1918}. This is the number of irreducible characters of $S_m$, and so certainly $m(S_m)$ is bounded by this number. Written as a function of $|S_m|=m!=n$, this becomes (of the order of)
\[ \frac{\e^{a\sqrt{f(n)}}}{b\cdot f(n)},\]
where $f(n)=\log n/\log \log n$. By removing $b$ from the denominator, we get a function that is definitely larger than $m(S_m)$ for sufficiently large $m$. Taking logarithms yields
\begin{align*} \log m(S_m)&\leq a\sqrt{f(n)}-\log f(n)
\\ &=a\sqrt{\log n/\log\log n}-\log(\log n/\log\log n)
\\ &\leq a\sqrt{\log n/\log\log n}-\log\log n.\end{align*}
Thus
\begin{align*}\log(\log m(S_m)+\log\log n)&\leq \log\(a\sqrt{\log n/\log\log n}\)
\\&=\log a+\half(\log\log n-\log\log\log n)
\\&<\half\log\log n\end{align*}
for sufficiently large $n$. This proves the first assertion.

\bigskip

Moving on to the groups of Lie type, let $G$ be a group of Lie type of the form in the table above, and let $m$ be the dimension of the natural module for $G$ (so that $m=d$ for $\SL_d(q)$, $m=2d+1$ for $\POmega_{2d+1}(q)$, and so on). If $G$ is untwisted, write $n=q^{m^2}$, and if $G$ is special unitary, write $n=(q^2)^{m^2}$. In all cases, $n>|G|$ since $n$ is equal to the total number of $m\times m$ matrices over $\F_q$ (or $\F_{q^2}$ in the twisted case).

Let us firstly consider the groups $G=\PSL_d(q)$, with $n=q^{d^2}$. By \cite[Theorem 327]{hardywright}, $\phi(a)\geq a^\delta$ for any $\delta<1$ and all sufficiently large $a$. Therefore, for all sufficiently large $n$,
\[ m(G)\geq \frac{\phi(q^d-1)}{d^2(q-1)}\geq \frac{(q^d-1)^\delta}{qd^2\log q}\approx \frac{q^{d\delta-1}}{\log n}.\]
(The middle inequality holds for all $q$ (even $q=2$) since $q\log q/(q-1)>1$ for all $q\geq 2$.) Taking logarithms gives
\[ \log m(G)\geq \log \(\frac{q^{d\delta-1}}{\log n}\)=(d\delta-1)\log q-\log\log n,\]
and thus
\begin{align*}\log\(\log m(G)+\log\log n\)&\geq \log\((d\delta-1)\log q\)
\\ &>\log d+\log\delta-1+\log\log q
\\ &>\half\log\log n+\log\delta-1,\end{align*}
since $\log\log n=2\log d+\log\log q$.
The same argument works for the other classical groups, completing the proof of our claim.

\section{Representation Growth of $p$-Groups}
\label{propgroups}

For powerful $p$-groups (i.e., groups $G$ for which $G'\leq G^p$ if $p$ is odd and $G'\leq G^4$ if $p=2$) one can get very good bounds on the number of conjugacy classes.

\begin{lem}[{{\cite[Lemma 4.7(ii)]{shalev1995}}}]\label{conjclassespowerful} If $G$ be a powerful finite $p$-group, then
\[ k(G)\geq (1-p^{-1})|G|^{1/d},\]
where $d=d(G)$ is the number of generators of $G$.
\end{lem}

Using this, it is very easy to give a lower bound for powerful pro-$p$ groups, and in fact a slightly larger class of pro-$p$ groups.

\begin{prop}\label{dgenprop} Let $G$ be a $d$-generator pro-$p$ with FAb, and suppose that $G$ has powerful finite images of arbitrarily large order (in particular, if $G$ is a powerful pro-$p$ group). For all powers $n$ of the prime $p$,
\[ R_n(G)\geq cn^{2/d},\]
where $c=(1-p^{-1})$.
\end{prop}
\begin{pf} Let $N$ be a normal subgroup such that $G/N$ is powerful of order $p^m$ where $m$ is even. (Since quotients of powerful groups are powerful, we can do this for all even $m$.) We have $k(G/N)\geq c|G/N|^{1/d}$; each of the irreducible representations of $G/N$ is of dimension less than $p^{m/2}$, and so
\[ R_{p^{m/2}}(G)\geq ap^{m/d};\]
writing $n=p^{m/2}$ we get $R_n(G)\geq cn^{2/d}$.
\end{pf}

If one wants a result on the numbers $r_n(G)$ rather than $R_n(G)$, then this is easy now.

\begin{cor} Let $G$ be a $d$-generator pro-$p$ with FAb, and suppose that $G$ has powerful finite images of arbitrarily large order (in particular, if $G$ is a powerful pro-$p$ group). For infinitely many powers $n$ of the prime $p$,
\[ r_n(G)\geq \frac{2cn^{2/d}}{\log_pn},\]
where $c=(1-p^{-1})$.
\end{cor}
This follows simply because there are at most $(\log_pn)/2$ degrees of irreducible representations of $G$ at most $n$.

A similar result can be obtained for some groups that are not powerful, like the Nottingham group.

\begin{prop} Let $p$ be an odd prime and let $G$ be the Nottingham group over $\F_q$, where $q$ is a power of $p$. For all powers $n$ of $p$ we have
\[ R_n(G)\geq cn^{2/3p},\]
where $c=c(q)$ depends only on $q$.
\end{prop}
\begin{pf} By \cite[Theorem 1.2]{jaikin2003}, for any normal subgroup $N$ of $G$, we have $k(G/N)\geq c|G/N|^{1/3p}$, where $c$ depends only on $q$. The method of proof of Proposition \ref{dgenprop} now gives the result.
\end{pf}

For all finitely generated pro-$p$ groups, until recently only a logarithmic bound for $R_n(G)$ was possible. However, in \cite{jaikin2009un}, Jaikin-Zapirain proved the following theorem.

\begin{thm}[Jaikin-Zapirain {{\cite{jaikin2009un}}}]\label{jaikinamazing} There is a constant $c>0$ such that, for all finite $p$-groups $G$, we have
\[ k(G)> c\log_p|G|\frac{\log_p\log_p|G|}{\log_p\log_p\log_p|G|}.\]
\end{thm}

Using this result, it is very easy to prove Theorem \ref{repgrowthpgroup}, via the same methods used for Proposition \ref{dgenprop}.

\bigskip

Andrei Jaikin-Zapirain has suggested the following slight improvement to Theorem \ref{repgrowthpgroup}, if one relaxes the condition that all sufficiently large $n$ satisfy the bound to just infinitely many $n$: in this case, one can get
\[ R_{p^n} \geq cn\log_p n\]
for infinitely many $n$ and some constant $c$ (independent of the group $G$). To see this, suppose that $G$ is a finitely generated pro-$p$ group. If $G$ is $p$-adic analytic, then $G$ contains a powerful subgroup of finite index, and hence the result follows by Proposition \ref{dgenprop} (for any $c>0$). If $G$ is not virtually powerful, then all dimension subgroups are distinct (see \cite[Theorem 11.5]{ddms}); write $a_n=|G:D_{2^n}(G)|$. For infinitely many $n$ we have $a_n/a_{n-1}\geq a_{n-1}/2$; let $H=G/D_{2^n}$ for some such $n$. Note that we have
\[ |D_{2^{n-1}}(H)|\geq |H|^{1/3}.\]
Using the proof of \cite[Claim 3.5]{jaikin2009un} (which states that if $G$ is a finite $p$-group with maximal powerful normal subgroup $P$, then $k(G/\Phi(P))\geq pm\log_pm/24$, where $m=d(P)$), one sees that if $P$ is a powerful normal subgroup containing $D_{2^{n-1}}(H)$ (which is elementary abelian as $D_{2^n}(H)=1$) then $d(P)\geq {b/3}$, where $|H|=p^b$, and so the claim is proved.

%
%

\section{Constructing Groups of Slow Representation Growth}
\label{alternatinggroups}

Theorem \ref{theorema} states that there is a divergent subsequence to the sequence $r_n(G)$, and Proposition \ref{propRnGprofinite} states that the sequence $R_n(G)$ strays above $\log n(\log \log n)^{1-\ep}$ infinitely often. In some sense therefore there is a `global' lower bound to the representation growth of a profinite group. In another sense, however, that we will consider in this section, there is not.

In \cite{segal2001}, Segal constructed finitely generated profinite groups whose finite images are iterated wreath products of finite simple groups. A different kind of group was constructed by Kassabov and Nikolov \cite{kassabovnikolov2006}, with finite images direct products of finite simple groups.

More specifically, let $\ms S$ be any infinite collection of finite simple groups, where each group may appear with a finite multiplicity. In \cite{kassabovnikolov2006} it was proved that, under suitable conditions for the multiplicities, there is a finitely generated profinite group whose finite images are \emph{exactly} the finite direct products of elements of $\ms S$. (One such suitable condition that we will use later is that all elements of $\ms S$ have multiplicity $1$.)

The Kassabov--Nikolov examples have representation growths that are reasonably easy to compute. Using alternating groups of varying degrees, Kassabov and Nikolov constructed profinite groups $G$ for which $R_n(G)$ is bounded between $n^b$ and $n^{b+\ep}$ for any $b>0$ and and $\ep>0$ (and $n$ sufficiently large), so that the abcissca of convergence of the zeta function is exactly $b$. (It could be that the representation growth of $G$ is, for example, $n^b\log n$.)

For functions $f$ that are supermultiplicative (i.e., $f(x)f(y)\geq f(xy)$), grow faster than any polynomial, and are below $n!$, it was also proved \cite[Theorem 1.8(a)]{kassabovnikolov2006} that there is a finitely generated profinite group $G$ such that $R_n(G)/f(n)\to 1$ as $n\to \infty$. One can achieve even faster growth by inserting (say) elementary abelian groups underneath the copies of the alternating groups. Thus there can be no bounds on the rate at which the sequence $R_n(G)$ may grow.

Let $\ms S=\{A_{n_1},A_{n_2},\dots\}$ be a collection of alternating groups, with $n_i<n_{i+1}$. Suppose that sequence $(n_i)$ grows very quickly; more precisely, suppose that $n_i>\prod_{j<i}(n_j!)$. This condition implies that the representations of degree at most $n_i-2$ are all representations of the product of the first $i-1$ elements of $\ms S$. Therefore, for $i\geq 2$,
\[ R_{n_i-2}(G)=\prod_{j<i} k(A_{n_j})\approx\prod_{j<i} p(n_j)/2,\]
where again $p(m)$ denotes the number of partitions of $m$. (The number of conjugacy classes of $A_m$ is approximately $p(m)/2$.) The reason for the $n_i-2$ is that $A_{n_i}$ has no representations of degree less than $n_i-1$, and exactly one of degree $n_i-1$, at least if $n_i\geq 7$.

Given any non-decreasing function $f:\N \to \N$, that tends to infinity, we may construct a finitely generated profinite group $G$ such that $R_n(G)< f(n)$ for infinitely many $n$. To see this, simply define $G$ to be as follows: let $n_1=7$, and choose $n_2$ such that $k(A_{n_1})> f(n_2)$ (and also $n_2>n_1!$). This ensures that $R_{n_2-2}(G)=k(A_{n_1})$. We repeat the process, choosing $n_3$ such that $n_3> (n_1!)(n_2!)$ and $f(n_3)>k(A_{n_1})k(A_{n_2})$, and so on.

This process produces a finitely generated group $G$ such that $R_{n_i-2}(G)<f(n_i-2)$ for all $i$. Thus it is not posssible to produce a global lower bound, proving Theorem \ref{repgrowthminimum}. The most sensible statement that one can make about lower bounds for representation growth is to require that the function $R_n(G)$ be greater than $f(n)$ \emph{infinitely many times}. In the example we constructed above, for any non-decreasing, divergent $f$, we can of course choose the $n_i$ so that $R_n(G)<f(n)$ for arbitrarily large intervals in $\N$. We therefore cannot, given a divergent non-decreasing function $f:\N\to \N$, even give lower bounds on `the proportion of $\N$' (e.g., using a measure like density) for which any finitely generated profinite group $G$ satisfies $R_n(G)>f(n)$, for example for $f(n)=c\log n(\log \log n)^{1-\ep}$ as in Proposition \ref{propRnGprofinite}.

Notice that, since $R_{n_i-2}(G)$ is the product of partition functions (roughly) for the groups above, we actually have that these groups $G$ satisfy $R_n(G)/f(n)>1$ for infinitely many $n$, where $f(n)$ is of the form $\e^{\alpha\sqrt{\log n}}$ for some $\alpha>0$, and so are a long way off the bound in Proposition \ref{propRnGprofinite}.

\bigskip

\bigskip

\noindent\textit{Acknowledgement}: Some of the work done here was completed during the \emph{Batsheva seminar on Representation Growth}, which took place between the 21st and 26th of June, 2009. I wish to thank the organizers for the opportunity to enter discussions with the various attendees. I would also like to thank Martin Kassabov for many interesting discussions during the Berlin summer school on finite simple groups and algebraic groups, which took place between the 31st of August and the 11th of September, 2009. In particular, he explained more fully the ideas behind \cite{kassabovnikolov2006}. Of course, I would like to thank Andrei Jaikin-Zapirain for suggesting the proof of Proposition \ref{propRnGprofinite}.

\bibliography{references}

\end{document}